\documentclass[12pt, letterpaper]{amsart}
\usepackage{amsmath, amssymb, latexsym, amsthm, mathrsfs}

\usepackage[top=2in, bottom=1.5in, left=1in, right=1in]{geometry}

\newcommand{\ed}{

\subsection*{Acknowledgments}
We thank Vaja Tarieledze for bringing relevant bibliography to our
attention and for vivid discussions. We thank Franklin Tall for comments
which helped to improve the presentation of the second part of the
paper. We also thank Assaf Rinot for an interesting discussion regarding Lemma~\ref{remains}.
We are indebted to Justin Moore for bringing Lemma \ref{vytjag}
to our attention, and to the referee for his time and effort.
The second named author was supported by the FWF grant I 1209-N25 and by the
Austrian Academy of Sciences, through the APART Program.
A part of this work was carried out while the first
named author was on Sabbatical leave at the Weizmann Institute of Science.
This author thanks his hosts for their kind hospitality.

\end{document}}

\newcommand{\nc}[1]{\newcommand{#1}}
\usepackage{graphicx}
\makeatletter
\newcommand*\nrotarrowconstructor[2]{%
  \mathrel{\m@th\sbox\z@{$ #1 $}%
    \raisebox{1.3\dp\z@}{%
      \makebox[\wd\z@][c]{%
        \reflectbox{\rotatebox[origin=cB]{90}{$ #2 $}}%
        \kern0.32\wd\z@%
      }}}%
}
\newcommand\nuparrow{\nrotarrowconstructor\uparrow\nrightarrow}

\makeatother
\nc{\set}[2]{\{\, #1 : #2 \,\}}
\nc{\agame}{\alpha_2^\mathrm{game}}
\nc{\qn}[1]{\par\smallskip {\sf *** #1 ***}\par\smallskip}
\nc{\Pa}[8]{\bibitem{#1} {#2}, \emph{#3}, {#4} \textbf{#5} ({#6}), {#7}--{#8}.}
\newcommand{\uhr}{\upharpoonright}
\nc{\tPa}[5]{\bibitem{#1} {#2}, \emph{#3}, {#4}, to appear.}
\nc{\sPa}[4]{\bibitem{#1} {#2}, \emph{#3}, {#4}, submitted.}
\nc{\Bc}[9]{\bibitem{#1} {#2}, \emph{#3}, in: \textbf{#4} (#5), #6 #7, #8--#9.}
\nc{\fD}{\mathfrak{D}}
\nc{\fX}{\mathfrak{X}}
\nc{\Onbd}{\rmO_{\mathrm{nbd}}} 
\nc{\Omnb}{\Omega_{\mathrm{nbd}}} 
\nc{\od}{\mathfrak{od}}
\nc{\Setting}[7]{\xymatrix@R=4pt@C=7pt{#1\ar@{-}[r]&#2\ar@{-}[r]&#3\\&#4\ar@{-}[u]\\
#5\ar@{-}[uu]\ar@{-}[r] & #6\ar@{-}[u]\ar@{-}[r] & #7\ar@{-}[uu]}}
\nc{\mx}[1]{\begin{matrix}#1\end{matrix}}
\nc{\pr}{\op{pr}}
\nc{\plim}{p\txt{-}\lim}
\nc{\osc}{\op{osc}}
\nc{\Bgp}{{\Z^\N}}
\nc{\Cgp}{{{\Z_2}^\N}}
\nc{\Cite}[1]{\textbf{[#1]}}
\nc{\Next}[1]{{#1^+}}
\nc{\Fr}{\mathit{F\!r}}
\nc{\intvl}[2]{{[#1(#2),\allowbreak #1(#2\!+\!1))}}
\nc{\Bdd}{\mathbf{B}}
\nc{\Cp}{\op{C}_\mathrm{p}}
\nc{\Ax}{\mathsf{Ax}}
\nc{\Dfin}{\mathfrak{D}_\mathrm{fin}}
\nc{\gp}{\mbox{-\textit{\tiny gp}}}
\nc{\grbl}{{\mbox{\textit{\tiny gp}}}}
\nc{\bbP}{\mathbb{P}}
\nc{\bbT}{\mathbb{T}}
\nc{\BOfat}{\B_{\Omega_{\mathrm{fat}}}}
\nc{\Bgood}{\B_{\mathrm{good}}}
\nc{\compactN}{\cl{\mathbb{N}}}
\nc{\blocks}[2]{\op{cl}_{#2}(#1)}
\nc{\blocksplus}[2]{\op{cl}^+_{#2}(#1)}
\nc{\arx}[1]{\texttt{http://arxiv.org/math/#1}}
\nc{\bq}{\begin{quote}}
\nc{\eq}{\end{quote}}
\nc{\cl}[1]{\overline{#1}}
\nc{\CH}{the Continuum Hypothesis}
\nc{\MA}{Martin's Axiom}
\nc{\Bfat}{\B_\mathrm{fat}}
\nc{\inv}{^{-1}}
\nc{\Cantor}{{\{0,1\}^\N}}
\nc{\bP}{\mathbf{P}}
\nc{\bof}{\op{\fb}}
\nc{\bofF}{\bof(\cF)}
\nc{\sr}[2]{{\txt{$#1$\\#2}}}
\nc{\nsr}[2]{#1}
\nc{\N}{\mathbb{N}}
\nc{\NN}{{\N^{\N}}}
\nc{\ZN}{{\Z^{\N}}}
\nc{\NNup}{{\N^{\uparrow\N}}}
\nc{\PN}{{P(\N)}}
\nc{\roth}{{[\N]^{\oo}}}
\nc{\Fin}{{[\N]^{<\oo}}}
\nc{\ici}{{[\N]^{(\oo,\oo)}}}
\nc{\Inc}{{\compactN^{\uparrow\N}}}
\nc{\powInc}[1]{{\big(\Inc\big)^{#1}}}
\nc{\powFin}[1]{{\big(\Fin\big)^{#1}}}
\nc{\powPN}[1]{{\big(\PN\big)^{#1}}}
\nc{\NcompactN}{{\compactN^\N}}
\nc{\seq}[1]{\{#1\}_{n\in\N}}
\nc{\sseq}[1]{\set{#1}{n\in\N}}
\nc{\Uarrow}{\smash{\big\uparrow}}
\nc{\LE}{\preccurlyeq}
\nc{\GE}{\succcurlyeq}
\nc{\op}{\operatorname}
\nc{\im}{\op{im}}
\nc{\Span}{\op{span}}
\nc{\maxfin}{\op{maxfin}}
\nc{\ran}{\op{range}}
\nc{\iso}{\cong}
\nc{\Madd}{{\M}^\star}
\nc{\cI}{\mathcal{I}}
\nc{\cJ}{\mathcal{J}}
\nc{\scrA}{\mathscr{A}}
\nc{\scrB}{\mathscr{B}}
\nc{\scrC}{\mathscr{C}}
\nc{\scrD}{\mathscr{D}}
\nc{\A}{\forall}
\nc{\B}{\mathrm{B}}
\nc{\cB}{\mathcal{B}}
\nc{\bB}{\mathbf{B}}
\nc{\Ga}{\Gamma}
\nc{\BG}{\B_\Ga}
\nc{\BL}{\B_\Lambda}
\nc{\BT}{\B_\Tau}
\nc{\BTstar}{\B_{\Tau^*}}
\nc{\BO}{\B_\Omega}
\nc{\DO}{\cD_\Omega}
\nc{\KO}{\cK_\Omega}
\nc{\CG}{C_\Ga}
\nc{\CL}{C_\Lambda}
\nc{\CT}{C_\Tau}
\nc{\CTstar}{C_{\Tau^*}}
\nc{\CO}{C_\Omega}
\nc{\COgp}{C_{\Omega^{\grbl}}}
\nc{\CLgp}{C_{\Lambda^{\grbl}}}
\nc{\BOgp}{\B_{\Omega}^{\grbl}}
\nc{\BLgp}{\B_{\Lambda^{\grbl}}}
\nc{\sfC}{\mathsf{C}}
\nc{\sfD}{\mathsf{D}}
\nc{\bD}{\mathbf{D}}
\nc{\Tau}{\mathrm{T}}
\nc{\cA}{\mathcal{A}}
\nc{\cK}{\mathcal{K}}
\nc{\cD}{\mathcal{D}}
\nc{\cF}{\mathcal{F}}
\nc{\cS}{\mathcal{S}}
\nc{\cG}{\mathcal{G}}
\nc{\cY}{\mathcal{Y}}
\nc{\J}{\mathcal{J}}
\nc{\cL}{\mathcal{L}}
\nc{\cM}{\mathcal{M}}
\nc{\cN}{\mathcal{N}}
\nc{\cO}{\mathcal{O}}
\nc{\rmO}{\mathrm{O}}
\nc{\cP}{\mathcal{P}}
\nc{\Q}{\mathbb{Q}}
\nc{\bbR}{\mathbb{R}}
\nc{\cU}{\mathcal{U}}
\nc{\Union}{\bigcup}
\nc{\cV}{\mathcal{V}}
\nc{\cW}{\mathcal{W}}
\nc{\Z}{{\mathbb Z}}
\nc{\Impl}{\Rightarrow}
\long\def\forget#1\forgotten{}
\nc{\ft}{\mathfrak{t}}
\nc{\fb}{\mathfrak{b}}
\nc{\fc}{\mathfrak{c}}
\nc{\fd}{\mathfrak{d}}
\nc{\fg}{\mathfrak{g}}
\nc{\oo}{\infty}
\nc{\fr}{\mathfrak{r}}
\nc{\fu}{\mathfrak{u}}
\nc{\fh}{\mathfrak{h}}
\nc{\fp}{\mathfrak{p}}
\nc{\fj}{\mathfrak{j}}
\nc{\fs}{\mathfrak{s}}
\nc{\w}{\omega}
\nc{\x}{\times}
\nc{\Iff}{\Leftrightarrow}
\nc\comp{^{\text{\tt c}}}
\nc{\nin}{\notin}
\nc{\cat}{\hat{\ }}
\nc{\sub}{\subseteq}
\nc{\spst}{\supseteq}
\nc{\sm}{\setminus}
\nc{\as}{\subseteq^*}
\nc{\rest}{\restriction}

\nc{\<}{\langle}

\nc{\E}{\exists}
\nc{\dom}{\op{dom}}
\nc{\cov}{\op{cov}}
\nc{\add}{\op{add}}
\nc{\cof}{\op{cof}}
\nc{\cf}{\op{cf}}
\nc{\non}{\op{non}}
\nc{\unif}{\op{non}}
\nc{\COV}{\op{COV}}
\nc{\ADD}{\op{ADD}}
\nc{\COF}{\op{COF}}
\nc{\NON}{\op{NON}}
\nc{\impl}{\to}
\nc{\Lp}{\mathcal{L_\p}}
\nc{\Wlog}{without loss of generality}
\newtheorem{thm}{Theorem}[section]
\nc{\bthm}{\begin{thm}} \nc{\ethm}{\end{thm}}
\newtheorem{prop}[thm]{Proposition}
\nc{\bprp}{\begin{prop}} \nc{\eprp}{\end{prop}}
\newtheorem{fact}[thm]{Fact}
\nc{\bfct}{\begin{fact}} \nc{\efct}{\end{fact}}
\newtheorem{prob}[thm]{Problem}
\nc{\bprb}{\begin{prob}} \nc{\eprb}{\end{prob}}
\newtheorem{lem}[thm]{Lemma}
\nc{\blem}{\begin{lem}} \nc{\elem}{\end{lem}}
\newtheorem{claim}[thm]{Claim}
\nc{\bclm}{\begin{claim}} \nc{\eclm}{\end{claim}}
\newtheorem{cor}[thm]{Corollary}
\nc{\bcor}{\begin{cor}} \nc{\ecor}{\end{cor}}
\newtheorem{conj}[thm]{Conjecture}
\nc{\bcnj}{\begin{conj}} \nc{\ecnj}{\end{conj}}
\theoremstyle{definition}
\newtheorem{defn}[thm]{Definition}
\nc{\bdfn}{\begin{defn}} \nc{\edfn}{\end{defn}}
\theoremstyle{remark}
\newtheorem{rem}[thm]{Remark}
\nc{\brem}{\begin{rem}} \nc{\erem}{\end{rem}}
\newtheorem{cnv}[thm]{Convention}
\nc{\bcnv}{\begin{cnv}} \nc{\ecnv}{\end{cnv}}
\newtheorem{exam}[thm]{Example}
\nc{\bexm}{\begin{exam}} \nc{\eexm}{\end{exam}}
\nc{\bpf}{\begin{proof}} \nc{\epf}{\end{proof}}
\nc{\be}{\begin{enumerate}}
\nc{\ee}{\end{enumerate}}
\nc{\bi}{\begin{itemize}}
\nc{\itm}{\item}
\nc{\ei}{\end{itemize}}
\nc{\Subsection}[1]{\goodbreak\subsection*{#1}}

\forget
\forgotten

\nc{\sone}{\mathsf{S}_1}
\nc{\sfin}{\mathsf{S}_\mathrm{fin}}
\nc{\ufin}{\mathsf{U}_\mathrm{fin}}
\nc{\Split}{\mathsf{Split}}

\nc{\gone}{\mathsf{G}_1}    \nc{\gfin}{\mathsf{G}_\mathrm{fin}}

\title[Sheaf amalgamations in topological groups]{Arhangel'ski\u{\i} sheaf amalgamations in\\
topological groups}
\author[Tsaban]{Boaz Tsaban}
\address[Tsaban]{Department of Mathematics, Bar-Ilan University, Ramat Gan 5290002, Israel and
Department of Mathematics, Weizmann Institute of Science, Rehovot 7610001, Israel}
\email{tsaban@math.biu.ac.il}
\urladdr{http://math.biu.ac.il/\~{}tsaban}

\author[Zdomskyy]{Lyubomyr Zdomskyy}
\address[Zdomskyy]{Kurt G\"odel Research Center for Mathematical Logic, University of Vienna,
W\"ahringer Str.\ 25, 1090 Vienna, Austria}
\email{lzdomsky@logic.univie.ac.at}
\urladdr{http://www.logic.univie.ac.at/~lzdomsky}

\keywords{%
Amalgamation of convergent sequences,
$\alpha_1$ space,
$\alpha_{1.5}$ space,
L-space.
}

\subjclass{%
26A03, 
03E75 
}

\begin{document}

\begin{abstract}
We consider amalgamation properties of convergent sequences in topological groups
and topological vector spaces.
The main result of this paper is that,
for arbitrary topological groups, Nyikos's property $\alpha_{1.5}$ is equivalent to
Arhangel'ski\u{\i}'s formally stronger property $\alpha_1$.
This result solves a problem of Shakhmatov (2002), and its proof uses a new
perturbation argument.
We also prove that there is a topological space $X$ such that the space
$\Cp(X)$ of continuous real-valued functions on $X$,
with the topology of pointwise convergence,
has Arhangel'ski\u{\i}'s property $\alpha_1$ but is not countably tight.
This result follows from results of Arhangel'ski\u{\i}--Pytkeev, Moore and Todor\v{c}evi\'c,
and provides a new solution, with stronger properties than the earlier solution, 
of a problem of Averbukh and Smolyanov (1968) concerning topological vector spaces.
\end{abstract}

\maketitle

\bigskip

\section{Sheaf amalgamations in topological groups}

To avoid trivialities, by \emph{convergent sequence} $x_n\to x$ we mean a proper one, that is,
such that $x\neq x_n$ for all $n$. This way, convergence is a property of countably infinite sets:
a countably infinite set $A$ converges to $x$ if all (equivalently, some) bijective enumerations of
$A$ converge to $x$. Thus, in the following definition,
by \emph{sequence} we always mean a countably infinite set.
The following concepts are due to Arhangel'ski\u{\i} \cite{Arhan72, Arhan81}, except
for $\alpha_{1.5}$ which is due to Nyikos \cite{Nyikos92}.

\bdfn\label{alphadef}
A topological space $X$ is $\alpha_i$, for $i=1,1.5,2,3,4$, if, respectively, for each $x\in X$ and
all pairwise disjoint sequences $S_1,S_2,\dots\sub X$, each converging to $x$,
there is a sequence $S\sub\Union_nS_n$ such that $S$ converges to $x$, and
\bi
\itm[($\alpha_1$)] $S_n\sm S$ is finite for all $n$.
\itm[($\alpha_{1.5}$)] $S_n\sm S$ is finite for infinitely many $n$.
\itm[($\alpha_2$)] $S_n\cap S$ is infinite for all $n$.
\itm[($\alpha_3$)] $S_n\cap S$ is infinite for infinitely many $n$.
\itm[($\alpha_4$)] $S_n\cap S$ is nonempty for infinitely many $n$.
\ei
\edfn

A survey of these properties is available in \cite{Shakhmatov02}.
In the integer-indexed properties $\alpha_i$, we may remove the requirement that the sequences
$S_1,S_2,\dots$ are pairwise disjoint \cite{Nyikos92}.
Indeed, we can move to subsequences $S_n'=S_n\sm\Union_{k<n}S_k$ of $S_n$, for $n\in\N$.
If $S_n'$ is infinite for infinitely many $n$, we can dispose
of the other ones.
And if not, then the sequence $S:=\Union_{k<n} S_k$, for any $n$ with $S_n'$ finite,
would be as required in $(\alpha_1)$.
However, removing the disjointness requirement in the property
$\alpha_{1.5}$ renders it superfluous: Applying it to the modified sequence
$\Union_{k\le n}S_k$, for $n\in\N$, the obtained sequence $S$ would be as required in $\alpha_1$.

Each of the properties in Definition \ref{alphadef} implies the subsequent one. To see that $\alpha_{1.5}$ implies $\alpha_2$,
for each $n$ decompose $S_n=\Union_k S_{nk}$, and take $S_n'=\Union_{m\le n} S_{mn}$ \cite{Nyikos92}.

None of the implications
$$\alpha_1\Impl \alpha_{1.5}\Impl \alpha_2\Impl \alpha_3\Impl \alpha_4$$
can be (provably) reversed. Not even in the class of Fr\'echet--Urysohn spaces \cite{Shakhmatov02}.
Recall that a topological space $X$ is \emph{Fr\'echet--Urysohn} if each point in the closure of a set
is in fact a limit of a sequence in that set.

In the present paper, we consider these properties in the context of \emph{topological groups}.
This direction was pioneered by Nyikos in his 1981 paper \cite{Nyikos81}.
In his paper, Nyikos proved that Fr\'echet--Urysohn groups are $\alpha_4$,
and that sequential $\alpha_2$ groups are Fr\'echet--Urysohn.
Shakhmatov \cite{Shakhmatov90} constructed, in the Cohen reals model, an example of a
Fr\'echet--Urysohn group which is not $\alpha_3$, and a Fr\'echet--Urysohn $\alpha_2$ group
which is not $\alpha_{1.5}$.
In particular, none of the implications
$$\alpha_1\Impl \alpha_2\Impl \alpha_3\Impl \alpha_4$$
is provably reversible in the realm of topological groups.
The question whether $\alpha_{1.5}$ groups are $\alpha_1$ is implicit
in Shakhmatov's paper.
The problem whether \emph{Fr\'echet--Urysohn} $\alpha_{1.5}$ groups are $\alpha_1$ is stated there.
This variant of the problem was settled in the positive by Shibakov, in his 1999 paper \cite{Shibakov99}.

In his 2002 chapter for \emph{Recent Progress in Topology} \cite{Shakhmatov02},
Shakhmatov cites Shiba\-kov's solution, and writes: ``It seems unclear if $\alpha_{1.5}$ and $\alpha_1$
are equivalent for all (i.e., not necessarily Fr\'echet--Urysohn) topological groups.''
For groups of the form $\Cp(X)$, the continuous real-valued functions on a space $X$, with the topology
of pointwise convergence, Sakai solved this problem in the positive \cite{SakaiRamsey}.
One step in his solution, uses a pullback method which was used earlier by Scheepers
\cite{Scheepers98} to show that for spaces of the form
$\Cp(X)$, we have that $\alpha_2=\alpha_3=\alpha_4$:
Replace the $n$-th sequence $\set{f_{nm}}{m\in\N}$ by
$\set{|f_{1m}|+\dots+|f_{nm}|}{m\in\N}$.
This approach is not applicable to arbitrary topological groups. Indeed, Sakai proves some of
his lemmata in the context of general topological groups, but his main theorems are proved
only in the case of $\Cp(X)$. The following theorem answers Shakhmatov's question.


\bthm\label{main}
A topological group is $\alpha_{1.5}$ if, and only if, it is $\alpha_1$.
\ethm
\bpf
Let $G$ be a topological group, and $S_1,S_2,\dots\sub G$ be sequences converging to $e$.
Let $T$ be any sequence converging to $e$ (e.g., let $T:=S_1$).
For each $n$, fix a bijective enumeration $S_n=\set{g_{nm}}{m\in\N}$.

Let $\set{(n_k,m_k)}{k\in\N}$ be an enumeration of the set $\N\x\N$ where each pair
$(n,m)$ appears infinitely often.
For each $k$, as the set $(T\sm\{t_1,\dots,t_{k-1}\})\cdot g_{n_km_k}$ is infinite, we can pick
an element
$$t_k\in T\sm\{t_1,\dots,t_{k-1}\}$$
such that
$$t_k\cdot g_{n_km_k}\notin\{t_1\cdot g_{n_1m_1},\ \dots\ ,\ t_{k-1}\cdot g_{n_{k-1}m_{k-1}}\}.$$
For each pair $(n,m)$, let $\set{k(n,m,i)}{i\in\N}$ be an increasing enumeration of
the set $\set{k}{(n_k,m_k)=(n,m)}$.
Note that the function $(n,m,i)\mapsto k(n,m,i)$ is injective.
For each $i$, define the following perturbation of $S_n$:
$$S_n^{(i)}=\{t_{k(n,1,i)}\cdot g_{n1}, t_{k(n,2,i)}\cdot g_{n2}, t_{k(n,3,i)}\cdot g_{n3}, \dots \}.$$
The sequence $S_n^{(i)}$ converges to $e$.
By the construction, the sets $S_n^{(i)}$, for $n,i\in\N$, are pairwise disjoint,
and therefore so are the sets
$$S_n'=S^{(n)}_1\cup S^{(n)}_2\cup \dots \cup S^{(n)}_n,$$
for $n\in\N$. Being finite unions of sequences converging to $e$, the sequences
$S_1',S_2',\dots$ converge to $e$, too.

Apply $\alpha_{1.5}$ to the sequences $S_1',S_2',\dots$, to find a sequence $S'$ converging to $e$,
such that the set $S_n'\sm S'$ is finite for each $n$ in an infinite set $I\sub\N$.
Define
$$S:=\Union_{n\in I}\Union_{j=1}^n \set{g_{jm}}{m\in\N, t_{k(j,m,n)}\cdot g_{jm}\in S'}.$$
Since for each $n\in I$ and each $j=1,\dots,n$,
we have that $t_{k(j,m,n)}\cdot g_{jm}\in S'$ for all but finitely many $m$,
the set $S_j\sm S$ is finite for all $j$.

Finally, note that $S$ is obtained by taking a subsequence of $S'$ and multiplying its elements by
distinct elements  $t_{k(j,m,n)}^{-1}$, that is elements of a subsequence of $\set{t^{-1}}{t\in T}$,
which also converges to $e$.
Thus, $S$ converges to $e$, too.
\epf

We obtain a short proof of a result of Nogura and Shakhmatov.

\bdfn[Nogura--Shakhmatov \cite{NoSh95}]
A topological space $X$ is \emph{Ramsey} if, whenever $\lim_n\lim_m x_{nm}=x$,
there is an infinite $I\sub\N$ such that for each neighborhood $U$ of $x$,
there is $k$ such that $\set{x_{nm}}{k<n<m,\ n,m\in I}\sub U$.
\edfn

In general $\alpha_1$ topological spaces need not be Ramsey.
In the context of topological groups, the above definition simplifies to the following one.

\blem[Sakai \cite{SakaiRamsey}]\label{RamGp}
A topological group $G$ is Ramsey if, and only if, whenever $\lim_m g_{nm}\allowbreak=e$ for all $n$,
there is an infinite set $I\sub\N$ such that the sequence
$\set{g_{nm}}{n,m\in I, n<m}$ converges to $e$.
\elem
\bpf
Assume that $\lim_m g_{nm}=g_n$ and $\lim_n g_n=e$.
For each $n$, define $g_{nm}'=g_n\inv g_{nm}$. Then $\lim_m g_{nm}=e$ for all $n$.
\epf

\bthm[Nogura--Shakhmatov \cite{NoSh95}]
Every $\alpha_{1.5}$ topological group is Ramsey.
\ethm
\bpf
Let $G$ be an $\alpha_{1.5}$ topological group. 
We establish the property stated in Lemma~\ref{RamGp}.
Assume that $\lim_m g_{nm}\allowbreak=e$ for all $n$.
By Theorem \ref{main}, $G$ is $\alpha_1$, and thus
there is an increasing function 
$f\colon\N\to\N$ such that the sequence $\set{g_{nm}}{m\ge f(n)}$ converges to $e$.
Take $I$ to be the image of $f$.
\epf

\section{New amalgamations}\label{sec:nam}

Using the above-mentioned pullback method of Scheepers,
Sakai proved that for groups of the form $\Cp(X)$,
\emph{Ramsey} is equivalent to $\alpha_2$ \cite{SakaiRamsey}.
The following problem, though, remains open.

\bprb[Shakhmatov \cite{Shakhmatov02}]\mbox{}
\be
\itm Is every (Fr\'echet--Urysohn) $\alpha_2$ topological group Ramsey?
\itm Is every (Fr\'echet--Urysohn) Ramsey topological group $\alpha_2$?
\ee
\eprb

In the forthcoming Definitions (\ref{dfn:1}, \ref{dfn:2}, \ref{dfn:3}, and \ref{dfn:4}),
we introduce several new local properties related to Ramsey and $\alpha_2$, 
and prove implications among them. The exact relations among these new
properties and among them and the classic ones remain unknown. Some 
of the most interesting problems that remain open are summarized in Section~\ref{sec:prbs}.

\bdfn\label{dfn:1}
A topological space $X$ is \emph{locally Ramsey} if,
for each $x\in X$, whenever $\lim_m x_{nm}=x$ for all $n$,
there is an infinite set $I\sub\N$ such that the sequence
$\set{x_{nm}}{n,m\in I, n<m}$ converges to $x$.
\edfn

Locally Ramsey spaces are $\alpha_3$.
By Lemma \ref{RamGp}, a topological group is Ramsey if, and only if, it is locally Ramsey.

\bdfn\label{dfn:2}
A topological space $X$ is $\alpha_{2^-}$ if, for each $x\in X$, whenever $\lim_m x_{nm}=x$ for all $n\in\N$,
there are natural numbers $m_1<m_2<\cdots$ such that the sequence
$\Union_n\{x_{1m_n},\dots,x_{nm_n}\}$ converges to $x$.
\edfn

Thus, every $\alpha_{2^-}$ topological space is $\alpha_2$.

\bprp
\mbox{}
\be
\itm Every $\alpha_{2^-}$ topological space is locally Ramsey.
\itm Every $\alpha_{2^-}$ topological group is Ramsey.
\ee
\eprp
\bpf
(1) Take $m_1<m_2<\cdots$ as in the definition of $\alpha_{2^-}$, and set $I:=\sseq{m_n}$.

(2) By (1) and Lemma \ref{RamGp}.
\epf

\bdfn\label{dfn:3}
A topological space $X$ is $\alpha_{3^-}$ if, for each $x\in X$, whenever $\lim_m x_{nm}=x$ for all $n$,
there are infinite sets $I,J\sub\N$ such that the sequence
$\set{x_{nm}}{n\in I, m\in J, n<m}$ converges to $x$.
\edfn

Thus, every locally Ramsey space is $\alpha_{3^-}$, and every $\alpha_{3^-}$ space is $\alpha_3$.
The above-mentioned results of Scheepers and Sakai follow.

\bprp
For topological groups of the form $\Cp(X)$, the properties
\begin{quote}
$\alpha_{2^-},\alpha_2,\alpha_{3^-},\alpha_3,\alpha_4$, locally Ramsey, and Ramsey,
\end{quote}
are equivalent.
\eprp
\bpf
By the above observations, it suffices to show that $\alpha_4$ implies $\alpha_{2^-}$ for
such spaces. This follows from Scheepers's pullback method:
Given sequences $S_n=\set{f_{nm}}{m\in\N}$, each converging to $0$,
replace each sequence $S_n$ with
$$S_n'=\set{|f_{1m}|+\dots+|f_{nm}|}{m\ge n}.$$
Applying $\alpha_4$ and thinning out, we obtain an increasing sequence
of indices $m_1<m_2<\cdots$ such that the sequence
$$|f_{1m_n}|+\dots+|f_{nm_n}|\quad (n\in\N)$$
converges to $0$.
Then the sequence $\Union_n\{f_{1m_n},\dots,f_{nm_n}\}$ converges to $0$.
\epf

\bdfn\label{dfn:4}
Let $X$ be a topological space, and $x\in X$.
The game $\agame(X,x)$ is played by two players, ONE and TWO,
and has an inning per each natural number. On the $n$th inning, ONE chooses
a sequence $S_n$ converging to $x$, and TWO responds by choosing a subsequence $T_n\sub S_n$.
TWO wins if the sequence $\Union_n T_n$ converges to $x$. Otherwise, ONE wins.
\edfn

\bprp
Assume that for each $x\in X$, ONE does not have a winning strategy in $\agame(X,x)$.
Then the space $X$ is $\alpha_{2^-}$ (and thus locally Ramsey).
\eprp
\bpf
Assume that $\lim_m x_{nm}=x$ for all $n$.
Consider the following strategy for ONE:
In the first inning, ONE plays the sequence $\set{x_{1m}}{m\in\N}$.
If TWO plays the subsequence
$$\set{x_{1m}}{m\in I_1},$$
then ONE responds by playing the sequence
$$\set{x_{2m}}{m\in I_1\sm\{\min I_1\}}.$$
In general, if in the $n$th inning TWO chooses a subsequence
$$\set{x_{nm}}{m\in I_n},$$
then ONE plays the sequence
$$\set{x_{n+1,m}}{m\in I_n\sm\{\min I_n\}}.$$
Since this strategy is not winning for ONE, there is a play lost
by ONE. Let $I_1,I_2,\dots$ be the infinite sets of sequence indices, which correspond
to the moves of TWO in this play. Define $m_n:=\min I_n$ for each $n$.
Then, for each $n$,
$$\Union_{n\in\N}\{x_{1m_n},\dots,x_{nm_n}\}\sub \Union_{n\in\N}T_n,$$
and the latter sequence converges to $x$.
\epf

\bcor
Let $G$ be a topological group.
If ONE does not have a winning strategy in the game $\agame(G,e)$,
then $G$ is Ramsey (indeed, $\alpha_{2^-}$).\qed
\ecor

\bprp
Let $X$ be an $\alpha_1$ space. For each $x\in X$,
ONE does not have a winning strategy in the game $\agame(X,x)$.
\eprp
\bpf
Define the game $\alpha_1^\mathrm{game}(X,x)$ corresponding to the property $\alpha_1$ (at $x$).
This game is similar to $\agame(X,x)$, with the only difference that here,
TWO must choose a \emph{cofinite} subset of each sequence provided by ONE.

\blem
A topological space $X$ is $\alpha_1$ if, and only if, for each point $x\in X$,
ONE does not have a winning strategy in the game $\alpha_1^\mathrm{game}(X,x)$.
\elem
\bpf
$(\Leftarrow)$ Immediate.

$(\Impl)$ The following method was used by Scheepers in \cite{coc1} to prove similar results
for games involving open covers.

Fix a strategy for ONE in $\alpha_1^\mathrm{game}(X,x)$.
For each sequence played by ONE, there are only countably many possible legal
responds by TWO. Let $\cF$ be the family of all possible sequences which ONE
may play according to the fixed strategy.
As the family $\cF$ is countable, we can apply $\alpha_1$ to $\cF$, and find for each
sequence $S\in\cF$ a
cofinite subset $S'\sub S$, such that the sequence $\Union_{S\in\cF}S'$ converges to $x$.

Consider a play where TWO responds to each given sequence $S_n$ by the sequence $S_n'$.
This play is lost by ONE.
\epf
If ONE does not have a winning strategy in $\alpha_1^\mathrm{game}(X,x)$,
then ONE does not have one in $\agame(X,x)$, where the moves of TWO
are less restricted.
\epf

\section{Sheaf amalgamations in topological vector spaces}\label{sec:tvs}

In their 1968 paper \cite{AvSmo68}, Averbukh and Smolyanov asked whether every
$\alpha_1$  topological vector space is Fr\'echet--Urysohn.
The problem was only settled in Plichko's 2008 paper  \cite{Plichko08}, using Banach spaces with
certain weak topologies. 
Knowledge that was available in the field of selection principles, even before its
solidification in 1996 \cite{coc1, coc2}, was enough to have a consistent counterexample for the
Averbukh--Smolyanov problem:
Assume \CH{}, and let $S\sub\bbR$ be a Sierpi\'nski set, that is,
a set of size continuum whose intersection with every Lebesgue-null set is countable.
It is known that every Borel image of a Sierpi\'nski set in $\NN$ is bounded, and consequently
the space $\Cp(S)$ is $\alpha_1$. On the other hand, the space $\Cp(S)$ cannot be Fr\'echet--Urysohn
since the set $S$ is not Lebesgue-null \cite{GN}.
Moreover, there is an example based solely on cardinality: It is known that the combinatorial
cardinal $\fp$ (respectively, $\fb$) is the minimal cardinality of a set $X\sub\bbR$
such that the space $\Cp(X)$ is not Fr\'echet--Urysohn (respectively, $\alpha_1$).
Thus, the consistent assumption $\fp<\fb$ provides a counterexample in a trivial manner.
We show in Theorem \ref{CpL} that this approach provides a counterexample, 
within ZFC. Moreover, this example has the following remarkable properties:
Every separable subspace is metrizable, but the topological vector space is not even
countably tight.

In the proof of Theorem \ref{CpL}, we will use several known facts, for which we provide
proofs for completeness.

\forget
For a product $\prod_{i\in I}X_i$ and $J\subset I$,
$\pr_J\colon\prod_{i\in I}X_i\to \prod_{i\in J}X_i$ is projection on the coordinates in $J$, that is,
restriction to $J$.

\blem[folklore]\label{ochev1}
Assume that $X$ is a hereditarily Lindel\"of subspace of a product space $\prod_{i\in I}X_i$,
$Y$ is a second countable Hausdorff space, and $f\colon X\to Y$ is continuous.
Then there are a countable set $J\sub I$ and a continuous function
$g\colon\prod_{i\in J}X_i\to Y$ such that $f=g\circ\pr_{J}$.
\elem
\begin{proof}
Let $\cB$ be a countable base for the topology of $Y$. For
each $B\in\cB$, find a countable family $\cU_B$ of
standard basic open sets in $\prod_{i\in I}X_i$ such that
$\Union\cU_B\cap X=f^{-1}(B)\cap X$. Let $J$ be the union of the supports of
all $U\in\cU_B$, $B\in\cB$.
The set $J$ is countable.
As $Y$ is Hausdorff, for all $x_0,x_1\in X$ with $f(x_0)\neq f(x_1)$, $\pr_J(x_0)\neq\pr_J(x_1)$.
Thus, we can define a function $g\colon\prod_{i\in J}X_i\to Y$ by $g(\pr_J(x))=f(x)$ for all $x\in X$.
By the choice of $J$, the function $g$ is continuous.
\end{proof}
\forgotten

General versions of the following fact were proved in the 1970's
(e.g., \cite{Kom78} and references therein).
Recall that the \emph{$\Sigma$-product} of spaces $X_i$, for $i\in I$,
with respect to a point $x\in \prod_{i\in I}X_i$, is the subspace
$\Sigma_{i\in I}X_i$ of the product space $\prod_{i\in I}X_i$,
consisting of all $y\in \prod_{i\in I}X_i$ such that
$y_i=x_i$ for all but countably many $i\in I$.

\bprp\label{ochev2}
Let $X$ be a $\Sigma$-product of a family of first countable spaces. Then:
\be
\itm Every countable subspace of $X$ is first countable.
\itm The space $X$ is $\alpha_1$.
\itm The space $X$ has countable tightness.
\itm The space $X$ is Fr\'echet--Urysohn.
\ee
\eprp
\bpf
(1) Countable subspaces of $X$ are supported on a countable set of indices.

(2) follows from (1).

(3) Let $X=\Sigma_{i\in I}X_i$, $A\sub X$ and $y\in\cl{A}$.
For each $i\in I$, let $\cB_i$ be a countable base at $y_i$.
For a finite set $F\sub I$ and an element $U\in\prod_{i\in F} \cB_i$, let
$$[U] := \set{x\in X}{\forall i\in F, x_i\in U_i}.$$

Fix an arbitrary, countably infinite set $I_1\sub I$.
Continue by induction on $n$.
Let $A_n\sub A$ be a countable set intersecting
$[U]$ for all finite $F\sub I_n$ and all $U\in \prod_{i\in F}\cB_i$.
Let $I_{n+1}$ be the union of $I_n$ and the supports of the elements of $A_n$.

The point $y$ is in the closure of the countable set $\Union_nA_n$.
Indeed, let $F$ be a finite subset of $I$, and $U\in \prod_{i\in F}\cB_i$.
Let $F_1=F\cap\Union_n I_n$, and $F_2=F\sm \Union_n I_n$.
As the set $F$ is finite, there is $n$ such that $F_1\sub I_n$.
Let $V=(\,U_i : i\in F_1\,)$.
Then there is an element $a\in A_n\cap[V]$. As the support of $a$ is contained in $I_{n+1}$,
$a_i=y_i$ for all $i\in F_2$.
Thus, $a\in [U]$.

(4) follows from (3) and (1).
\epf


The following result, brought to our attention by J. Moore,
is proved for S in \cite[Theorem~7.10]{Tobook},
where it is pointed out that the L case is analogous.
For completeness, we provide a proof for the L case, which is the one
needed here.

\blem\label{vytjag}
Assume that $Y$ is a regular topological space with all finite powers Lindel\"of and countably tight,
and $X$ is a non-separable subspace of $Y$.
There exists a c.c.c.\ poset $\mathbb P$ such that, in $V^{\mathbb P}$,
the space $X$ has an uncountable discrete subspace.
\elem
\bpf
It suffices to show that there are a c.c.c.\ poset $\mathbb P$ and a family $\mathcal D$ of $\aleph_1$ many
dense subsets of $\mathbb P$ such that:
\begin{quote}
For each ZFC model $V'\spst V$ with $\omega_1^{V'}=\omega_1^V$,
if there is in $V'$ a filter $G\sub\mathbb P$ meeting each $D\in\mathcal D$,
then the space $X$ has an uncountable discrete subspace in $V'$.
\end{quote}
Passing to a subset of $X$, if necessary, we may assume that
$X=\set{x_\xi}{\xi<\omega_1}$ and
$\overline{\set{x_\xi}{\xi<\alpha}}^Y\cap\set{x_\eta}{\eta\geq\alpha}=\emptyset$
for every $\alpha<\omega_1$. There are two cases to consider.

\subsubsection*{Case 1}
$\overline{\set{x_\xi}{\xi<\alpha}}^Y\cap\overline{\set{x_\eta}{\eta\geq\alpha}}^Y=\emptyset$
for all $\alpha<\omega_1$; in other words, $X$ is a free sequence
in $Y$. Since the space $Y$ has countable tightness,
$\overline{X}^{Y}=\bigcup_{\alpha<\omega_1}\overline{\set{x_\beta}{\beta<\alpha}}^{Y}$.
The space $\overline{X}^{Y}$ is closed in $Y$, and thus Lindel\"of.
On the other hand, the family
$$\set{\overline{X}^{Y}\setminus \overline{\set{x_\beta}{\beta\geq\alpha}}^{Y}}{\alpha<\omega_1} $$
is an open cover of $\overline{X}^{Y}$ without a countable subcover; a contradiction.

\subsubsection*{Case 2}
$\overline{\set{x_\xi}{\xi<\alpha}}^Y\cap\overline{\set{x_\eta}{\eta\geq\alpha}}^Y\neq\emptyset$
for some $\alpha$. In particular, the set $\set{x_\eta}{\eta\geq\alpha}$ is not
compact. We may assume that the space $X$ is not compact.
Let $\cU$ be an ultrafilter on $X$ whose elements
are uncountable. If there exists some $\alpha$ such that $\cU$ contains all open neighborhoods in $X$ of $x_\alpha$, then the
Hausdorff property implies that every  $x_\beta$ for
$\beta\neq\alpha$ has a neighborhood in $X$ which is not in
$\cU$. By removing a point from $X$, if needed,
we may assume that every element of $X$ has a neighborhood  in $X$
that is not in $\cU$.

For each $\alpha$, pick neighborhoods $U_\alpha,V_\alpha$ of $x_\alpha$ in $Y$ such that
$\overline{V_\alpha}\sub U_\alpha$,  $\overline{U_\alpha}\cap
\overline{\set{x_\xi}{\xi<\alpha}}^Y=\emptyset$, and
$\set{U_\alpha\cap X}{\alpha<\omega_1}\sub P(X)\setminus\cU$.
Then finitely many sets $U_\alpha$ cannot cover a co-countable subset of $X$.
Let $\mathbb P$ be the poset consisting of all
finite sets $\{\alpha_0,\ldots,\alpha_{n-1}\}\sub\omega_1$, $\alpha_0<\cdots<\alpha_{n-1}$,
such that $x_{\alpha_j}\notin V_{\alpha_i}$ whenever $i<j$.
A condition $H$ is stronger than $F$, $H\leq F$, if $F\sub H$.

Assume, towards a contradiction, that there is an uncountable antichain
$\set{F_\alpha}{\alpha<\omega_1}$ in $\mathbb P$.
For incompatible elements $F,H\in\mathbb P$, the elements
$F\setminus H$ and $H\setminus F$ are also incompatible.
By the $\Delta$-System Lemma, we may assume that
the sets $F_\alpha$ are pairwise disjoint,
$\min F_\alpha>\max F_\beta$ for all $\beta<\alpha$, and $|F_\alpha|=n$
for all $\alpha$. For each $\alpha$,
let $\{\xi^0_\alpha,\ldots,\xi^{n-1}_\alpha\}$ be
the increasing enumeration of $F_\alpha$.
Set
\begin{eqnarray*}
W^0_\alpha:=\set{(x_0,\ldots,x_{n-1})\in X^{n}}{\forall i,j<n, (x_i\notin U_{\xi^j_\alpha})},\\
W^1_\alpha:=\set{(x_0,\ldots,x_{n-1})\in X^{n}}{\exists i,j<n, (x_i\in\overline{V_{\xi^j_\alpha}})}.
\end{eqnarray*}
Then $W^0_\alpha\cap W^1_\alpha=\emptyset$, and the sets $W^0_\alpha,W^1_\alpha$ are closed.
Moreover, by our choice of the sets $U_\delta$, we have that
$(x_{\xi^i_\beta})_{i<n}\in W^0_\alpha$ for all $\beta<\alpha$.
By the definition of $\mathbb P$ and the incompatibility of the sets $F_\alpha$ and $F_\beta$,
we have that $(x_{\xi^i_\beta})_{i<n}\in W^1_\alpha$ for all $\beta>\alpha$.
Thus, the subset
$A:=\set{\vec{x}_\alpha:=(x_{\xi^i_\alpha})_{i<n}}{\alpha<\omega_1}$ of
$X^n$ satisfies
$$\overline{\set{\vec{x}_\beta}{\beta<\alpha}}^{Y^n}\cap
\overline{\set{\vec{x}_\beta}{\beta\geq\alpha}}^{Y^n}=\emptyset$$
for all $\alpha$; a contradiction.

Thus, the forcing notion $\mathbb P$ is c.c.c.
For each $\alpha<\omega_1$,
let $D_\alpha:=\set{F\in\mathbb{P}}{\max F>\alpha}$.
Since no finite subfamily of
$\set{U_\alpha}{\alpha<\omega_1}$ covers a co-countable subset of $X$,
each set $D_\alpha$ is dense in $\mathbb P$. Assume that $G$ is a
subfilter of $\mathbb P$ (possibly, in some extension $V'\spst V$)
which intersects every set $D_\alpha$.
Then $x_\beta\notin V_{\alpha}$ for all $\beta,\alpha\in\Union G$: if
$\beta<\alpha$ this follows from the choice of $V_\alpha$, and if
$\beta>\alpha$ this follows from the existence of an element
$F\in G$ containing both $\alpha$ and $\beta$.
Thus, $G$ gives rise to the  discrete subspace $\set{x_\alpha}{\alpha\in\Union G}$  of
$X$, which is uncountable if $\omega_1^{V'}=\omega_1^V$.
\epf

We are ready for the main result of this section.
An \emph{L-space} is a hereditarily Lindel\"of nonseparable
topological space. The existence of L-spaces was established by
Moore in \cite{Moore06}.
A classical result of Arhangel'ski\u{\i} and, independently, Pytkeev,
asserts that a function space $\Cp(X)$ has countable tightness if and only if all
finite powers of the space $X$ are Lindel\"of.

\bthm\label{CpL}
There is a hereditarily Lindel\"of nonseparable Fr\'echet--Urysohn space $L$, such that:
\be
\item The space $\Cp(L)$ is $\alpha_1$;
moreover, every separable subspace of $\Cp(L)$ is metrizable.
\item The space $\Cp(L)$ is not Fr\'echet--Urysohn;
moreover, it is not countably tight.
\ee
\ethm
\bpf
Let $L$ be an L-space of the kind constructed by Moore \cite{Moore06}.
Following Todor\v{c}evi\'c \cite{Todorcevic89},
Moore considered a function $\osc\colon\set{(\alpha,\beta)\in\w_1^2}{\alpha<\beta}\to\w$
with strong combinatorial properties.
Let $(z_\alpha)_{\alpha<\w_1}$ be a sequence of rationally independent
points on the multiplicative
circle group $\bbT=\set{z\in\mathbb C}{|z|=1}$.
For each $\beta<\w_1$, define an element $w_\beta\in\bbT^{\w_1}$ by
$$w_\beta(\alpha)=
\begin{cases}
z_\alpha^{\osc(\alpha,\beta)+1} & \alpha<\beta \\
1 & \mbox{otherwise}
\end{cases}
$$
By Theorem 7.11 of \cite{Moore06}, the set $L=\set{w_\beta}{\beta<\w_1}$ is an L-space.
By Theorem 7.8 of \cite{Moore06}, or directly by Proposition \ref{ochev2},
the space $L$ is Fr\'echet--Urysohn.

(1)  Let $D$ be a countable subset of $\Cp(L)$.
Since $L$ is a hereditarily Lindel\"of subspace of a product space and $\bbR$ is
a second countable Hausdorff space, 
every continuous function $f\colon L\to \bbR$ is determined by
countably many coordinates; equivalently, 
there are $\alpha<\w_1$ and a continuous function $g_\alpha\colon\pr_\alpha[L]\to\bbR$
such that $f=g_\alpha\circ \pr_\alpha$.

\blem
For each $\alpha<\w_1$, the set $\pr_\alpha[L]$ is countable.
\elem
\bpf
By \cite[Proposition~7.13]{Moore06}, the subtree
$\set{\osc(\cdot,\delta)\uhr\alpha}{\delta\geq\alpha}$ of
$\omega^{<\omega_1}$ is Aronszajn, where
$\osc(\cdot,\delta)\colon\xi\mapsto \osc(\xi,\delta)$ for
$\xi<\delta$. By the definition of Aronszajn tree, the set
$$\set{\osc(\cdot,\delta)\uhr\alpha}{\alpha<\delta<\omega_1}$$
is countable for each $\alpha<\omega_1$.
Thus, the set $\set{w_\delta\uhr\alpha}{\alpha<\delta<\omega_1}$ is countable,
and hence so is the set $\pr_\alpha[L]$.
\epf

As the set $D$ is countable, there is $\alpha<\w_1$ such that
every function $f\in D$ is determined by a continuous function on
the first $\alpha$ coordinates.
Thus, the function
\begin{eqnarray*}
\pr_\alpha^* \colon \Cp(\pr_\alpha[L]) & \to & \Cp(L)\\
 g & \mapsto & g\circ\pr_\alpha
\end{eqnarray*}
is an embedding (e.g., \cite[Proposition~0.4.6]{Arhan92}).
As the set $\pr_\alpha[L]$ is countable, the space $\Cp(\pr_\alpha[L])$ is metrizable,
and therefore so is its image, which contains $D$.

(2) By Lemma \ref{ochev2}, every finite power of $\Sigma_{\alpha<\w_1}\bbT$ is countably tight.
As countable tightness is hereditary, all finite powers of $L$ are countably tight.
By Lemma~\ref{vytjag} with $X=Y=L$, we have that if all
finite powers of $L$ are Lindel\"of, then  there is a c.c.c.\ poset
$\mathbb P$  such that  $L$ has an uncountable discrete subspace in
$V^{\mathbb P}$. But in the proof of \cite[Theorem~7.17]{Moore06},
it is pointed out that the space $L$ remains an L-space in c.c.c.\ forcing
extensions. In fact, c.c.c.\ is not necessary, as the following lemma shows.

\blem\label{remains}
Moore's L-space remains an L-space in every forcing extension that does not
collapse $\aleph_1$.
\elem
\bpf[Proof sketch]
In accordance with \cite[Definition~2.1]{Moore06},
the construction of $L$ is based on a $C$-sequence
$$\bar{C}=\langle\, C_\alpha:\alpha<\omega_1,\alpha\mbox{ limit}\,\rangle.$$
The function $\osc$ is constructed from $\bar{C}$ in a way
that, for each poset $\mathbb P$ preserving $\omega_1$, the
constructions of $\osc$ in $V$ and in $V^{\mathbb P}$ give
the same function, and hence give rise to the same subspace of the
$\Sigma$ product of circles.\footnote{This can be checked by going through the relevant definitions in \cite[\S{}~4]{Moore06},
without involving any deep absoluteness arguments.}
By the same proof  carried out in
$V^{\mathbb P}$, this space is an $L$ space in $V^{\mathbb P}$.
\epf

It follows that some finite power of $L$ is not Lindel\"of. By the Arhangel'ski\u{\i}--Pytkeev Theorem,
the space $\Cp(X)$ is not countably tight.
\epf

\section{Open problems and closing remarks}\label{sec:prbs}

By Section~\ref{sec:nam} we have that, for topological groups,
$$\alpha_1\Iff\alpha_{1.5}\Impl \mathrm{ONE}\nuparrow\agame(G,e)
\Impl \alpha_{2^-}\Impl \mbox{Ramsey}\Impl\alpha_{3^-}\Impl\alpha_3\Impl\alpha_4$$
and $\alpha_{2^-}\Impl \alpha_2\Impl\alpha_3$.

\bprb
Are there, in ZFC or consistently, topological groups $G$ that are
\be
\itm $\alpha_{3}$ but not $\alpha_{3^-}$?
\itm $\alpha_{3^-}$ but not Ramsey?
\itm Ramsey but not $\alpha_{2^-}$?
\itm $\alpha_{2^-}$ but ONE has a winning strategy in $\agame(G,e)$?
\itm not $\alpha_1$ and ONE has no winning strategy in $\agame(G,e)$?
\itm $\alpha_2$ but not $\alpha_{2^-}$?
\itm $\alpha_{3^-}$ but not $\alpha_2$?
\ee
\eprb

\bprb
Let $X$ be a Tychonoff space such that the function space
$\Cp(X)$ satisfies $\alpha_2$.
Does it follow that ONE does not have a winning strategy in the game $\agame(\Cp(X),0)$?
\eprb

The results and methods of Section~\ref{sec:tvs} are already used in a number of papers,
including \cite{AD, AT, Peng, ST, Tall}.
The direct union of an $L$-space and the Sorgenfrey line is an $L$-space with non-Lindel\"of square.
However, such a space does not enjoy the properties described in Section~\ref{sec:tvs}.
Answering a question from an earlier version of this paper, Yinhe Peng proved that
the square of Moore's original L-space is also non-Lindel\"of \cite{Peng}.
Peng's arguments are highly nontrivial, and use the fine details of Moore's construction.
Our proof in item 2 of Theorem \ref{CpL} is potentially more general, as it applies to
all absolute modifications of Moore's construction where Lemma \ref{remains} holds.
We do not know any modification of Moore's construction where Lemma~\ref{remains} fails.

\ed